\documentclass[10pt,reqno]{amsart}
\usepackage{amsmath}
\usepackage{amsthm}
\usepackage{amssymb}
\usepackage{amsfonts}
\usepackage{latexsym}
\usepackage{hyperref}
\usepackage{cite}

\newcommand{\norm}[1]{\| #1 \|}

\newcommand{\N}{\mathbb{N}}

\renewcommand{\phi}{\varphi}

\theoremstyle{plain}

\theoremstyle{definition}

\allowdisplaybreaks

\begin{document}
\bibliographystyle{plain}

    \title{On the norm closure problem for complex symmetric operators}

    \author{Stephan Ramon Garcia}
    \address{   Department of Mathematics\\
            Pomona College\\
            Claremont, California\\
            91711 \\ USA}
    \email{Stephan.Garcia@pomona.edu}
    \urladdr{http://pages.pomona.edu/\textasciitilde sg064747}

	\author{Daniel E.~Poore}


    \thanks{Partially supported by National Science Foundation Grant DMS-1001614.}

    \begin{abstract}
    	We prove that the set of all complex symmetric operators on a separable, infinite-dimensional
	Hilbert space is not norm closed.
    \end{abstract}

\maketitle

In \cite[Sect.~3]{CSPI}, it is asked 
whether the set of all complex symmetric
operators on a separable, infinite-dimensional Hilbert space is norm closed.
We answer this question in the negative.
Let $S(a_0,a_1,a_2,\ldots) = (0,a_0,a_1,\ldots)$ denote the unilateral shift on $\ell^2(\N)$ and let $\cong$
denote unitary equivalence.  Note that
\begin{equation*}
	T_n
	\quad =\quad \tfrac{n}{n+1} S \oplus (\bigoplus_{ \substack{j=1 \\ j \neq n} }^{\infty} \tfrac{j}{j+1} S )
	\oplus (\bigoplus_{j=1}^{\infty} \tfrac{j}{j+1} S^* )
	\quad \cong\quad \bigoplus_{j=1}^{\infty} \tfrac{j}{j+1} (S \oplus S^*)
\end{equation*}
is complex symmetric by \cite[Ex.~5]{CSO2}.
On the other hand, $T_n$ converges in norm to 
\begin{equation*}
	T\quad=\quad S \oplus (\bigoplus_{j=1 }^{\infty} \tfrac{j}{j+1} S)\oplus (\bigoplus_{j=1 }^{\infty} \tfrac{j}{j+1} S^*)
	\quad \cong\quad S \oplus \bigoplus_{j=1}^{\infty} \tfrac{j}{j+1} (S \oplus S^*).
\end{equation*}
Since $\norm{S^k(1,0,0,\ldots)} = 1$, there
is an $x$ so that $\norm{T^k x} = 1$ for $k\geq 0$.
However,
\begin{equation*}
T^*= S^* \oplus \bigoplus_{j=1}^{\infty} \tfrac{j}{j+1} (S^* \oplus S) = S^* \oplus \text{(a strict contraction)}
\end{equation*}
possesses no such vector since $(S^*)^k$ tends strongly to zero.  This precludes the existence
of a conjugation $C$ (i.e., an isometric, conjugate-linear involution)
such that $T = CT^*C$.  Thus $T$ is not complex symmetric. \qed

\medskip	
\noindent\textbf{Remarks}:  
We thank D.~Sherman for his helpful suggestions.  We also note that S.~Zhu, C.G.~Li, and Y.Q.~Ji
discovered a different approach \cite{ZLJ} shortly before us.

\bibliography{ONCPCSO}

\end{document}